\documentclass[11pt]{article}
\usepackage[T1]{fontenc}
\usepackage[utf8]{inputenc}
\usepackage{lmodern}
\usepackage[letterpaper,margin=1.1in]{geometry}
\usepackage{microtype}
\usepackage{amsmath,amssymb,amsthm,mathtools}
\usepackage{graphicx}
\usepackage{booktabs,array}
\usepackage[font=small,labelfont=bf,labelsep=period,skip=6pt]{caption}
\usepackage{float}
\usepackage[section]{placeins}
\usepackage[table]{xcolor}
\usepackage{enumitem}
\usepackage{tikz}
\usetikzlibrary{arrows.meta,calc,positioning,fit,backgrounds}
\usepackage[colorlinks=true,linkcolor=accent,citecolor=accent,urlcolor=accent]{hyperref}
\usepackage[capitalise,noabbrev]{cleveref}

\definecolor{accent}{HTML}{1C5CAB}
\definecolor{cblue}{HTML}{2A78D6}
\definecolor{corange}{HTML}{EB6834}
\definecolor{ink}{HTML}{0B0B0B}
\definecolor{muted}{HTML}{52514E}
\definecolor{rule}{HTML}{C9C8C2}
\definecolor{caqua}{HTML}{1BAF7A}
\definecolor{nucfill}{HTML}{EEF4FC}
\definecolor{cytfill}{HTML}{FBFAF7}
\definecolor{bloodfill}{HTML}{FDEEE8}

\newcommand{\E}[1]{\left\langle #1\right\rangle}
\newcommand{\dd}{\mathrm{d}}
\newcommand{\dv}[2]{\frac{\dd #1}{\dd #2}}
\newcommand{\hill}[1]{\left(\frac{K}{K+#1/V}\right)}
\newcommand{\spn}{\sigma_{pn}}
\newcommand{\ph}{\,\mathrm{h}^{-1}}
\newcommand{\D}{\mathcal{D}}
\DeclareMathOperator{\re}{Re}
\newtheorem{proposition}{Proposition}
\theoremstyle{remark}
\newtheorem{remark}{Remark}

\numberwithin{equation}{section}
\crefname{equation}{equation}{equations}

\tikzset{
  >={Stealth[length=5pt,width=4.2pt]},
  flow/.style={->,line width=0.8pt,draw=ink},
  cat/.style={->,line width=0.8pt,draw=muted,dash pattern=on 3pt off 2pt},
  inhib/.style={-{Bar[width=7pt,line width=0.9pt]},line width=0.8pt,draw=corange},
  sp/.style={circle,draw=ink,line width=0.7pt,fill=white,inner sep=0pt,minimum size=22pt,font=\small},
  rate/.style={font=\footnotesize,text=ink,inner sep=1.5pt},
  lab/.style={font=\footnotesize\sffamily,text=muted},
  sink/.style={font=\small,text=muted,inner sep=1pt},
}
\newcommand{\cellpic}[2]{%
  \node[sp] (#2mn) at (0,0.75) {$m_n#1$};
  \node[sp] (#2pn) at (0,-0.75) {$p_n#1$};
  \node[sp] (#2mc) at (3.0,0.75) {$m_c#1$};
  \node[sp] (#2pc) at (3.0,-0.75) {$p_c#1$};
  \begin{scope}[on background layer]
    \node[draw=muted,line width=0.6pt,fill=cytfill,rounded corners=9pt,inner xsep=14pt,inner ysep=16pt,
          fit={(-1.95,1.3) (#2mc) (#2pc) (0,-0.75)}] (#2cell) {};
    \node[draw=cblue,line width=0.6pt,fill=nucfill,rounded corners=7pt,inner xsep=9pt,inner ysep=9pt,
          fit={(-1.85,0.75) (#2mn) (#2pn)}] (#2nuc) {};
  \end{scope}
  \node[lab,anchor=south west] at ($(#2nuc.north west)+(0.02,0.0)$) {nucleus};
  \node[lab,anchor=south east] at ($(#2cell.north east)+(-0.15,-0.42)$) {cytoplasm};
  \node[font=\scriptsize\sffamily,text=muted] (#2dna) at (-1.55,0.75) {DNA};
  \draw[flow] (#2dna) -- node[rate,above] {$\alpha P_0#1$} (#2mn);
  \draw[flow] (#2mn) -- node[rate,above] {$\gamma_m$} (#2mc);
  \draw[cat] (#2mc) -- node[rate,right] {$\beta$} (#2pc);
  \draw[flow] (#2pc) -- node[rate,below] {$\gamma_p$} (#2pn);
  \node[sink] (#2sinkp) at (-1.55,-0.75) {$\varnothing$};
  \draw[flow] (#2pn) -- node[rate,below] {$\delta_p$} (#2sinkp);
  \draw[inhib] (#2pn) to[bend left=28] node[rate,left,pos=0.42,text=corange] {$\xi,\eta$} ($(#2dna)!0.55!(#2mn)+(0,-0.06)$);
}

\title{\bfseries Spontaneous Synchronization in Cellular Circadian Clocks}
\author{Yuelong Li\\[2pt]\small Courant Institute of Mathematical Sciences,\\\small New York University, New York, NY 10012, USA}
\date{December 2022}
\begin{document}
\maketitle

\begin{abstract}
\noindent We study a population of identical cellular circadian oscillators, each a four-species transcription--translation feedback loop of the type analysed by Wang and Peskin, coupled by the exchange of cytoplasmic mRNA with a shared, well-mixed bloodstream at rate $\epsilon$. Simulations started from uniformly random phases show a collective damping of the mean rhythm followed by spontaneous synchronization. We show that the blood compartment tracks the population mean on the time scale $1/[N_c(\epsilon+\delta_m)]$, and that the population means obey the single-cell equations with a correction proportional to the variance $\spn^2$ of nuclear protein across cells; the correction is accurate only near synchrony. Deviations from the mean obey a linear system with periodic coefficients, whose leading Floquet exponent $-\kappa(\epsilon)$ is the asymptotic synchronization rate. Using the fact that the four first-order stages of the feedback loop commute, we prove that $\kappa(\epsilon)=\epsilon/4+O(\epsilon^2)$ on the limit cycle in the large-population limit. Numerically computed exponents agree with the decay rates measured in simulations to within $0.2\%$ for $5$ to $400$ cells. The approach to synchrony from an incoherent state is slower than $\kappa$, and at the weakest coupling tested the population remains in a two-cluster state. With intrinsic molecular noise, uncoupled cells lose their autocorrelation within about ten cycles, whereas coupled cells retain it and sustain a Kuramoto order parameter of $0.68$ to $0.93$.
\end{abstract}


\section{Introduction}

Beyond the one-directional flow of transcription and translation, proteins translated in the cytoplasm can be imported into the nucleus and repress their own gene. This delayed negative feedback is the core of cellular circadian clocks and goes back to the Goodwin oscillator~\cite{Goodwin}. Wang and Peskin~\cite{Wang} put forth a four-species model of this type, with mRNA and protein in the nucleus and in the cytoplasm, and studied its entrainment by light in the presence of molecular noise.

A single cell is a poor clock: its period drifts with molecular noise. Tissue-level rhythms rely on many cells keeping a common phase, and coupling through a shared diffusible signal is a standard mechanism for this. Populations of oscillators coupled through their mean field synchronize spontaneously~\cite{Winfree,Kuramoto,Strogatz}, and Gonze \emph{et al.}~\cite{Gonze} showed this for Goodwin-type circadian oscillators coupled through a common neurotransmitter concentration.

Here we extend the Wang--Peskin oscillator to $N_c$ identical cells that exchange cytoplasmic mRNA with a shared bloodstream. The coupling is built so that a synchronized population reproduces the single-cell dynamics exactly. We ask three questions: whether the population synchronizes from random initial phases, how the population average behaves on the way, and how fast synchrony is approached as a function of the exchange rate $\epsilon$.

The contributions are as follows.
\begin{itemize}[leftmargin=1.4em,itemsep=1pt,topsep=2pt]
\item A coupled model (\cref{sec:model}) and simulations showing collective damping followed by synchronization (\cref{sec:results}).
\item A reduction to the population means: the blood compartment is slaved to $\E{m_c}$, and the means obey the single-cell equations with a variance correction whose range of validity we measure (\cref{sec:mean}).
\item The asymptotic synchronization rate as a Floquet exponent, with a proof that it equals $\epsilon/4$ to first order, and a numerical comparison over exchange rates and population sizes (\cref{sec:rate}).
\item Stochastic simulations showing how coupling preserves the autocorrelation of individual cells under molecular noise (\cref{sec:noise}).
\end{itemize}

\section{Model}\label{sec:model}

\subsection{Single cell}

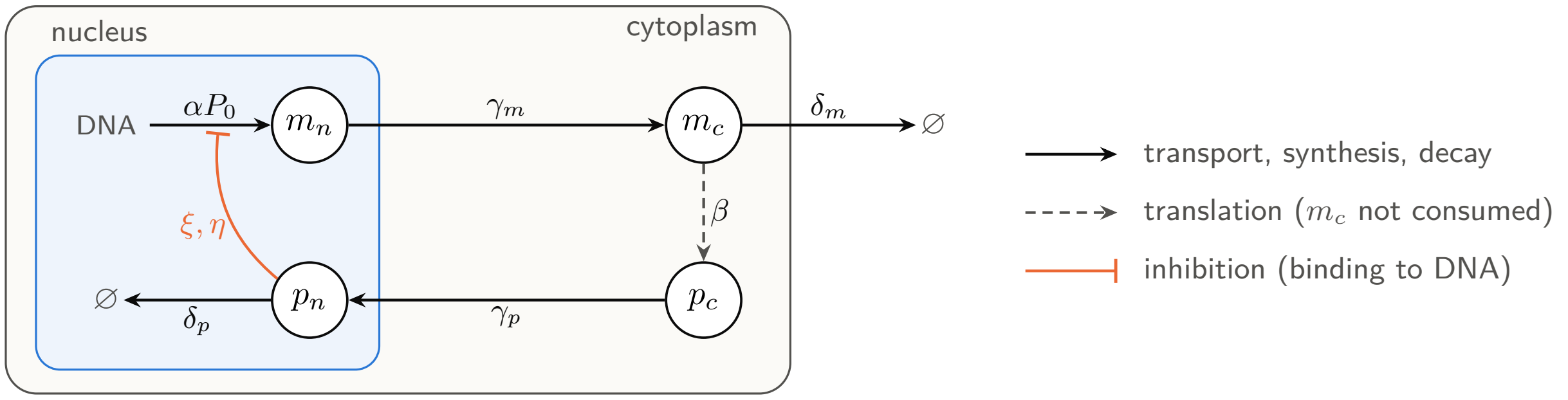
\begin{figure}[H]
\centering
\begin{tikzpicture}[x=1.35cm,y=1.2cm]
  \cellpic{}{a}
  \node[sink] (asinkm) at (4.75,0.75) {$\varnothing$};
  \draw[flow] (amc) -- node[rate,above] {$\delta_m$} (asinkm);
  \begin{scope}[shift={(5.45,0.5)}]
    \draw[flow] (0,0) -- (0.7,0); \node[lab,anchor=west] at (0.8,0) {transport, synthesis, decay};
    \draw[cat] (0,-0.5) -- (0.7,-0.5); \node[lab,anchor=west] at (0.8,-0.5) {translation ($m_c$ not consumed)};
    \draw[inhib] (0,-1.0) -- (0.7,-1.0); \node[lab,anchor=west] at (0.8,-1.0) {inhibition (binding to DNA)};
  \end{scope}
\end{tikzpicture}
\caption{The single-cell oscillator. Nuclear protein $p_n$ binds the DNA and shuts off transcription, closing a delayed negative-feedback loop.}
\label{fig:single-diagram}
\end{figure}

Each cell contains four species: mRNA in the nucleus ($m_n$), mRNA in the cytoplasm ($m_c$), protein in the cytoplasm ($p_c$) and protein in the nucleus ($p_n$). mRNA is transcribed in the nucleus and exported at rate $\gamma_m$; protein is translated from cytoplasmic mRNA at rate $\beta$ and imported into the nucleus at rate $\gamma_p$; cytoplasmic mRNA and nuclear protein are degraded at rates $\delta_m$ and $\delta_p$. Nuclear protein represses transcription by binding to $r$ sites on the DNA, with binding rate $\xi$ and unbinding rate $\eta$ per site (\cref{fig:single-diagram,tab:single-reactions}). All rates are probabilities per hour, and we work with molecule numbers rather than concentrations.

\begin{table}[H]
\centering\small
\caption{Reactions of the single-cell model. $\#(x)$ is the number of molecules of species $x$.}
\label{tab:single-reactions}
\begin{tabular}{@{}clll@{}}
\toprule
 & \textbf{Reaction} & \textbf{Rate} & \textbf{Result}\\
\midrule
1 & Transcription & $\alpha P_0$ & $\#(m_n)\rightarrow\#(m_n)+1$\\
2 & Export & $\gamma_m m_n$ & $\#(m_n)\rightarrow\#(m_n)-1$,\quad $\#(m_c)\rightarrow\#(m_c)+1$\\
3 & mRNA degradation & $\delta_m m_c$ & $\#(m_c)\rightarrow\#(m_c)-1$\\
4 & Translation & $\beta m_c$ & $\#(p_c)\rightarrow\#(p_c)+1$\\
5 & Import & $\gamma_p p_c$ & $\#(p_c)\rightarrow\#(p_c)-1$,\quad $\#(p_n)\rightarrow\#(p_n)+1$\\
6 & Protein degradation & $\delta_p p_n$ & $\#(p_n)\rightarrow\#(p_n)-1$\\
\bottomrule
\end{tabular}
\end{table}

Transcription proceeds only when no site is occupied, at rate $\alpha P_0$, where $P_l$ is the probability that $l$ sites are occupied. Binding is fast, so the states $l$ and $l+1$ are in equilibrium:
\begin{equation}
\xi\,(r-l)(p_n-l)\,P_l = \eta\,(l+1)\,P_{l+1}.
\end{equation}
With $K=k_2/k_1$ the ratio of the macroscopic unbinding and binding constants and $V$ the nuclear volume, $\eta/\xi=KV$. For $p_n\gg r$, normalizing $\sum_l P_l=1$ gives
\begin{equation}
P_0 = \left[\sum_{l=0}^r \binom{r}{l}\left(\frac{p_n}{KV}\right)^{l}\right]^{-1}=\hill{p_n}^r .
\label{eq:P0}
\end{equation}

In the continuous limit the molecule numbers obey
\begin{equation}
\begin{aligned}
\dv{m_n}{t} &= \alpha\hill{p_n}^r - \gamma_m m_n, &\qquad \dv{p_c}{t} &= \beta m_c - \gamma_p p_c,\\
\dv{m_c}{t} &= \gamma_m m_n-\delta_m m_c, & \dv{p_n}{t} &= \gamma_p p_c -\delta_p p_n .
\end{aligned}
\label{eq:single}
\end{equation}
Wang and Peskin~\cite{Wang} compare this limit with the stochastic model in detail. We write $h(p)=\alpha\,[K/(K+p/V)]^r$ for the transcription rate.

\subsection{Fixed point and onset of oscillation}\label{sec:hopf}

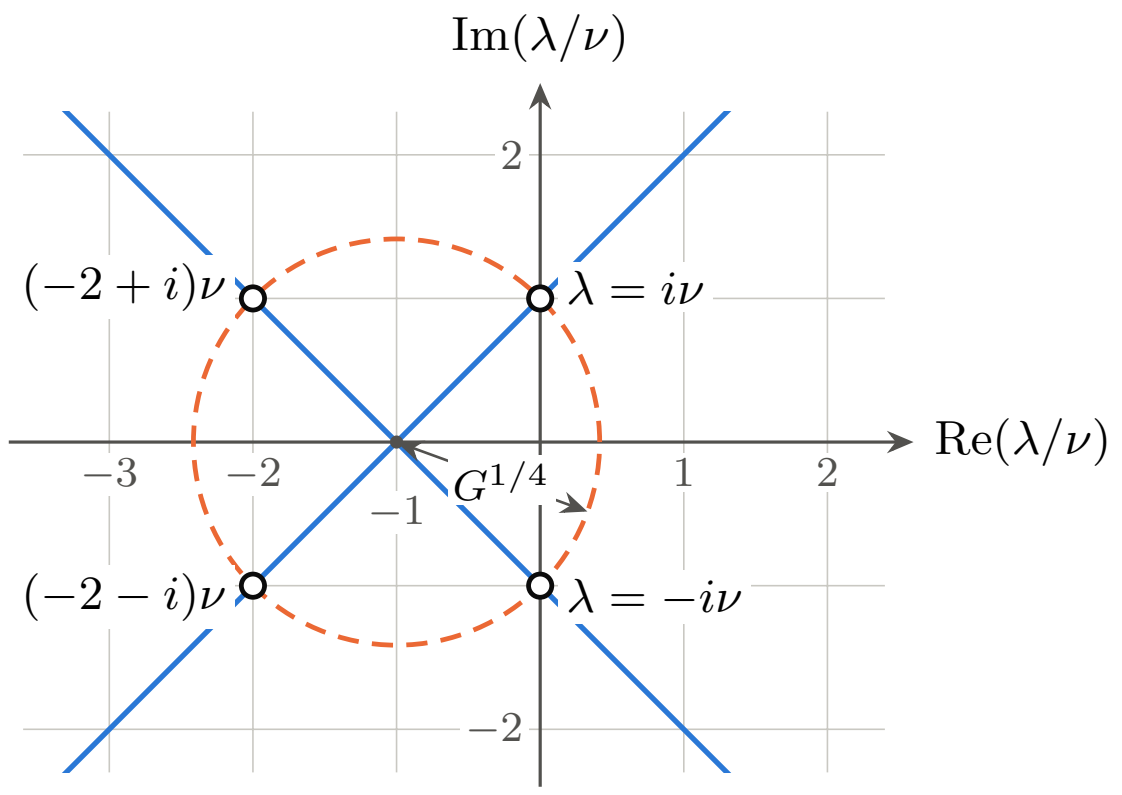
\begin{figure}[H]
\centering
\begin{tikzpicture}[scale=0.95]
  \draw[rule,line width=0.3pt,step=1] (-3.6,-2.3) grid (2.4,2.3);
  \draw[->,line width=0.6pt,draw=muted] (-3.7,0) -- (2.6,0) node[right,font=\footnotesize] {$\re(\lambda/\nu)$};
  \draw[->,line width=0.6pt,draw=muted] (0,-2.4) -- (0,2.5) node[above,font=\footnotesize] {$\operatorname{Im}(\lambda/\nu)$};
  \foreach \x in {-3,-2,1,2} \node[below,font=\scriptsize,text=muted,fill=white,inner sep=1pt,yshift=-2pt] at (\x,0) {$\x$};
  \node[below,font=\scriptsize,text=muted,fill=white,inner sep=1pt,yshift=-9pt] at (-1,0) {$-1$};
  \foreach \y in {-2,2} \node[left,font=\scriptsize,text=muted,fill=white,inner sep=1pt,xshift=-2pt] at (0,\y) {$\y$};
  \begin{scope}
    \clip (-3.6,-2.3) rectangle (2.4,2.3);
    \draw[line width=1pt,cblue] (-4,-3) -- (2,3);
    \draw[line width=1pt,cblue] (-4,3) -- (2,-3);
  \end{scope}
  \draw[line width=0.9pt,dash pattern=on 4pt off 2.5pt,corange] (-1,0) circle ({sqrt(2)});
  \foreach \p in {(0,1),(0,-1),(-2,1),(-2,-1)} \draw[line width=0.9pt,ink,fill=white] \p circle (2.3pt);
  \fill[muted] (-1,0) circle (1.3pt);
  \draw[<->,line width=0.5pt,draw=muted] (-1,0) -- node[font=\scriptsize,fill=white,inner sep=1pt,pos=0.55] {$G^{1/4}$} ({-1+sqrt(2)*cos(-20)},{sqrt(2)*sin(-20)});
  \node[font=\footnotesize,fill=white,inner sep=1.5pt,anchor=west] at (0.12,1.08) {$\lambda=i\nu$};
  \node[font=\footnotesize,fill=white,inner sep=1.5pt,anchor=west] at (0.12,-1.08) {$\lambda=-i\nu$};
  \node[font=\footnotesize,fill=white,inner sep=1.5pt,anchor=east] at (-2.12,1.08) {$(-2+i)\nu$};
  \node[font=\footnotesize,fill=white,inner sep=1.5pt,anchor=east] at (-2.12,-1.08) {$(-2-i)\nu$};
\end{tikzpicture}
\caption{Roots of the characteristic equation \eqref{eq:char} in the plane of $\lambda/\nu$. They lie on the lines at $45^\circ$ through $-1$ (solid), at distance $G^{1/4}$ from $-1$ (dashed circle). Drawn at the Hopf threshold $G=4$, where the leading pair is on the imaginary axis; for $G>4$ the circle grows and the leading pair moves into the right half-plane.}
\label{fig:eigen}
\end{figure}

The fixed point $p_{n0}$ of \eqref{eq:single} satisfies
\begin{equation}
\alpha\beta\hill{p_{n0}}^r = \delta_m\delta_p\, p_{n0},
\label{eq:equil}
\end{equation}
and the slope of the transcription rate there is
\begin{equation}
h'(p_{n0}) = -\frac{r\,\delta_p\delta_m}{\beta}\left(\frac{p_{n0}/V}{K+p_{n0}/V}\right)\equiv -a .
\label{eq:a}
\end{equation}
We take all four first-order rates equal, $\gamma_m = \delta_m = \gamma_p = \delta_p = \nu$. Linearizing \eqref{eq:single} about the fixed point then gives the characteristic equation
\begin{equation}
\left(\frac{\lambda}{\nu}+1\right)^4 = -\frac{a\beta}{\nu^2} = -G,
\qquad G\equiv r\left(\frac{p_{n0}/V}{K+p_{n0}/V}\right),
\label{eq:char}
\end{equation}
whose roots are $\lambda/\nu=-1+G^{1/4}e^{\pm i\pi/4}$ and $-1+G^{1/4}e^{\pm 3i\pi/4}$ (\cref{fig:eigen}). The leading pair has real part $\nu\,(G^{1/4}/\sqrt{2}-1)$, which is positive exactly when
\begin{equation}
G>4 .
\label{eq:hopf}
\end{equation}
At $G=4$ the leading pair is $\lambda=\pm i\nu$ and the fixed point undergoes a Hopf bifurcation; for $G>4$ it is unstable and trajectories settle on a limit cycle. Since $G<r$, oscillation requires $r\geq5$.

\subsection{Parameters}

\Cref{tab:params} lists the values used throughout. We fix the fixed point at $p_{n0}=500$, which keeps $p_n\gg r$ along the orbit as \cref{eq:P0} requires, and set $\alpha$ from \cref{eq:equil}:
\begin{equation}
\alpha = \frac{\nu^2\,p_{n0}}{\beta}\left(\frac{K+p_{n0}/V}{K}\right)^r.
\label{eq:alpha}
\end{equation}
The translation rate $\beta=10\ph$ follows~\cite{Wang}. With $K=63.5$\,pL$^{-1}$ the gain is $G=4.70$, above the threshold \eqref{eq:hopf}, and the leading eigenvalues at the fixed point are $\lambda=(0.041\pm1.041\,i)\,\nu$.

\begin{table}[H]
\centering\small
\caption{Parameter values.}
\label{tab:params}
\begin{tabular}{@{}llll@{}}
\toprule
\textbf{Symbol} & \textbf{Meaning} & \textbf{Value} & \textbf{Unit}\\
\midrule
$\nu$ & common rate $\gamma_m=\delta_m=\gamma_p=\delta_p$ & $2\pi/22 \approx 0.2856$ & h$^{-1}$\\
$\beta$ & translation rate & $10$ & h$^{-1}$\\
$r$ & number of DNA binding sites & $5$ & --\\
$V$ & nuclear volume & $0.5$ & pL\\
$K$ & dissociation constant $k_2/k_1$ & $63.5$ & pL$^{-1}$\\
$p_{n0}$ & fixed point of $p_n$ & $500$ & molecules\\
$\alpha$ & maximal transcription rate, \cref{eq:alpha} & $5.374\times10^{6}$ & h$^{-1}$\\
$G$ & loop gain, \cref{eq:char} & $4.70$ & --\\
\bottomrule
\end{tabular}
\end{table}

With $\alpha$ set by \cref{eq:alpha}, rescaling time by $\nu$ and the mRNA counts by $\nu^{-1}$ removes $\nu$ from \eqref{eq:single}, so the period is inversely proportional to $\nu$. For the values of \cref{tab:params} the limit cycle has $T\nu=6.667$, that is $T=23.34$\,h (\cref{fig:single}); a period of exactly 24\,h corresponds to $\nu=0.2778\ph$.

\subsection{Coupled cells}

\begin{figure}[H]
\centering
\begin{tikzpicture}[x=1.25cm,y=1.2cm]
  \begin{scope}[shift={(0,1.9)}]\cellpic{^1}{a}\end{scope}
  \begin{scope}[shift={(0,-2.0)}]\cellpic{^2}{b}\end{scope}
  \node[lab,anchor=east] at ($(acell.west)+(-0.1,0)$) {cell 1};
  \node[lab,anchor=east] at ($(bcell.west)+(-0.1,0)$) {cell 2};
  \node[font=\large,text=muted] at (1.0,-3.75) {$\vdots$};
  \node[lab,anchor=east] at ($(bcell.west)+(-0.1,-1.75)$) {cell $N_c$};
  \begin{scope}[on background layer]
    \fill[bloodfill,rounded corners=9pt] (5.0,-4.1) rectangle (6.9,3.75);
    \draw[corange,line width=0.6pt,rounded corners=9pt] (5.0,-4.1) rectangle (6.9,3.75);
  \end{scope}
  \node[lab,anchor=south] at (5.95,3.77) {bloodstream};
  \node[sp] (ms) at (5.95,0) {$m_s$};
  \foreach \c/\side in {a/north west, b/south west} {
    \coordinate (A) at (\c mc.east); \coordinate (B) at (ms.\side);
    \draw[flow] ($(A)!2.6pt!90:(B)$) -- node[rate,above,sloped] {$\mu_c$} ($(B)!2.6pt!-90:(A)$);
    \draw[flow,<-] ($(A)!2.6pt!-90:(B)$) -- node[rate,below,sloped] {$\mu_s$} ($(B)!2.6pt!90:(A)$);
  }
  \node[sink] (sinks) at (8.0,0) {$\varnothing$};
  \draw[flow] (ms) -- node[rate,above,pos=0.6] {$N_c\delta_m$} (sinks);
\end{tikzpicture}
\caption{The coupled model. Each cell exchanges cytoplasmic mRNA with the shared blood pool $m_s$ (rates $\mu_c$ out, $\mu_s$ in), where it is degraded.}
\label{fig:multi-diagram}
\end{figure}
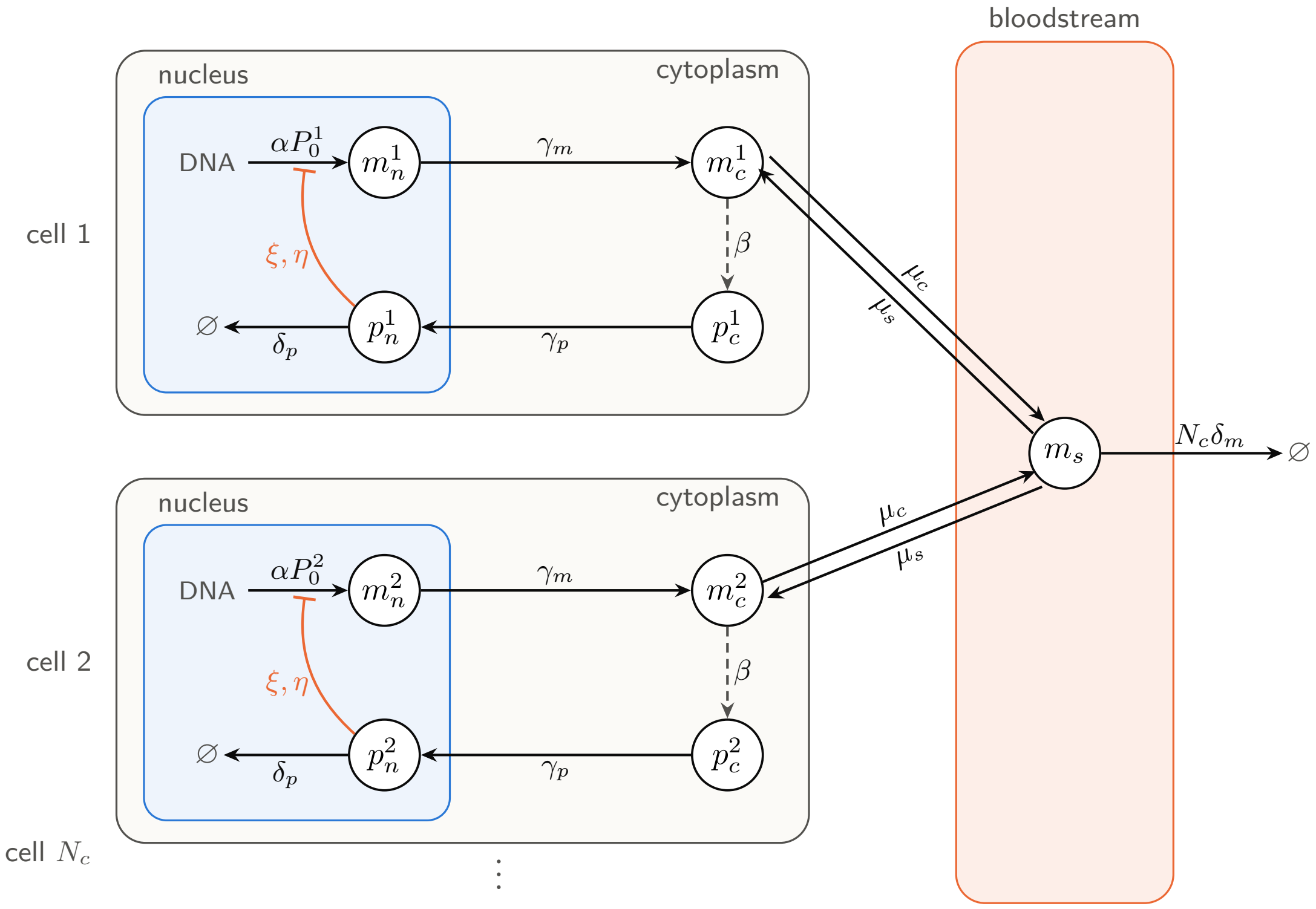

We attach $N_c$ identical cells to one bloodstream (\cref{fig:multi-diagram}). Cytoplasmic mRNA leaves each cell into the blood at rate $\mu_c$ and enters each cell from the blood at rate $\mu_s$; the blood pool $m_s$ is well mixed, which is reasonable because blood circulates much faster than the hours over which mRNA levels change. The cells thus form a star graph around the blood compartment. mRNA is no longer degraded inside the cells but in the blood, at rate $N_c\delta_m m_s$ (\cref{tab:multi-reactions}).

\begin{table}[H]
\centering\small
\caption{Reactions of the coupled model for cell $i=1,\dots,N_c$. Reactions that differ from \cref{tab:single-reactions} are in bold.}
\label{tab:multi-reactions}
\begin{tabular}{@{}clll@{}}
\toprule
 & \textbf{Reaction} & \textbf{Rate} & \textbf{Result}\\
\midrule
1 & Transcription & $\alpha P_0^i$ & $\#(m_n^i)\rightarrow\#(m_n^i)+1$\\
2 & Export & $\gamma_m m_n^i$ & $\#(m_n^i)\rightarrow\#(m_n^i)-1$,\quad $\#(m_c^i)\rightarrow\#(m_c^i)+1$\\
3 & \textbf{Diffusion} & $\mu_c m_c^i$ & $\#(m_c^i)\rightarrow\#(m_c^i)-1$,\quad $\#(m_s)\rightarrow\#(m_s)+1$\\
4 & Translation & $\beta m_c^i$ & $\#(p_c^i)\rightarrow\#(p_c^i)+1$\\
5 & Import & $\gamma_p p_c^i$ & $\#(p_c^i)\rightarrow\#(p_c^i)-1$,\quad $\#(p_n^i)\rightarrow\#(p_n^i)+1$\\
6 & Protein degradation & $\delta_p p_n^i$ & $\#(p_n^i)\rightarrow\#(p_n^i)-1$\\
7 & \textbf{Infusion} & $\mu_s m_s$ & $\#(m_c^i)\rightarrow\#(m_c^i)+1$,\quad $\#(m_s)\rightarrow\#(m_s)-1$\\
8 & \textbf{Blood degradation} & $N_c \delta_m m_s$ & $\#(m_s)\rightarrow\#(m_s)-1$\\
\bottomrule
\end{tabular}
\end{table}

The continuous equations are
\begin{equation}
\begin{aligned}
\dv{m_n^i}{t} &= h(p_n^i) - \gamma_m m_n^i, &\qquad \dv{p_c^i}{t} &= \beta m_c^i - \gamma_p p_c^i,\\
\dv{m_c^i}{t} &= \gamma_m m_n^i-\mu_c m_c^i+\mu_s m_s, & \dv{p_n^i}{t} &= \gamma_p p_c^i -\delta_p p_n^i,\\
\dv{m_s}{t} &= \mathrlap{\mu_c\sum_{i=1}^{N_c} m_c^i-N_c(\mu_s+\delta_m)\, m_s .}
\end{aligned}
\label{eq:multi}
\end{equation}
We introduce the exchange rate $\epsilon$ and set
\begin{equation}
\mu_c = \delta_m+\epsilon, \qquad \mu_s = \epsilon .
\label{eq:mu}
\end{equation}
With this choice a synchronized population in which $m_s=m_c^i$ has the extra export $\epsilon m_c$ cancelled by the import $\epsilon m_s$, leaving the single-cell degradation $\delta_m m_c$: synchronized cells follow \eqref{eq:single} and keep its period.

\section{Numerical methods}\label{sec:methods}

\paragraph{Integration.} \Cref{eq:single,eq:multi} are integrated with \texttt{scipy.integrate.solve\_ivp}~\cite{SciPy}, using the embedded Dormand--Prince RK5(4) pair~\cite{DormandPrince} with adaptive steps and relative and absolute tolerances of $10^{-10}$. As a check, the numerical derivative of the total protein $p_n+p_c$ along a single-cell solution agrees with $\beta m_c-\delta_p p_n$ to within $8\times10^{-4}\ph$, against rates of up to $289\ph$.

\paragraph{Initial conditions.} A single cell is started from $(m_n,m_c,p_c,p_n)=(10.9,\,1.88,\,250,\,500)$ and integrated until it has settled on the limit cycle. One cycle $\mathbf{S}(\phi)$, $\phi\in[0,2\pi)$, is cut between two consecutive minima of $m_n$ after $t=240$\,h (\cref{fig:single}). Each cell of a population starts at $\mathbf{S}(\phi_i)$ with $\phi_i$ drawn independently and uniformly from $[0,2\pi)$, and the blood starts empty, $m_s(0)=0$. A population therefore starts incoherent.

\paragraph{Observables.} We measure synchrony by the standard deviation of nuclear protein across cells,
\begin{equation}
\spn(t)=\E{\big(p_n^i-\E{p_n}\big)^2}^{1/2},\qquad \E{x}=\frac{1}{N_c}\sum_{i=1}^{N_c}x^i ,
\end{equation}
and by the Kuramoto order parameter $R(t)=\big|\E{e^{i\theta_i}}\big|$~\cite{Kuramoto}. The phase $\theta_i$ of a cell is the polar angle of $(p_n^i-p_{n0},\,p_c^i-p_{n0})$, mapped through the limit cycle so that it advances uniformly in time on it. $R=1$ for identical phases and $R\sim N_c^{-1/2}$ for random ones. For a time series $x(t)$ the normalized autocorrelation is
\begin{equation}
C(\tau)=\frac{\overline{\big(x(t)-\bar{x}\big)\big(x(t+\tau)-\bar{x}\big)}}{\overline{\big(x(t)-\bar{x}\big)^2}},
\label{eq:acf}
\end{equation}
with the bar a time average.

\paragraph{Floquet exponents.} The limit cycle is located to a tolerance of $10^{-11}$ by event detection on $\dd m_n/\dd t=0$. The variational equation \eqref{eq:variational} is integrated over one period from the identity matrix, and the exponents are $T^{-1}\log|\rho_k|$ for the eigenvalues $\rho_k$ of the resulting monodromy matrix. For finite $N_c$ the same is done around the synchronized orbit of one cell together with the blood compartment.

\paragraph{Decay rates.} Asymptotic decay rates are least-squares slopes of the cycle-averaged $\log\spn$ while it lies between $10^{-2}$ and $10^{-6}$. The transient rate is $\ln 10$ divided by the time $\spn$ takes to fall to a tenth of its initial value, averaged over five sets of initial phases.

\paragraph{Stochastic simulation.} The reactions of \cref{tab:multi-reactions} are simulated with integer molecule numbers by tau-leaping~\cite{Gillespie} with a step of $0.002$\,h. Production events (transcription, translation) are Poisson distributed, and first-order losses are binomial, so that counts cannot become negative. Halving the step leaves the time-averaged order parameter unchanged within its fluctuations ($0.72$ against $0.71$).

All runs are listed in \cref{app:setups}.

\section{Synchronization in simulations}\label{sec:results}

\subsection{Single cell}

\begin{figure}[!htbp]
\centering
\includegraphics[width=\textwidth]{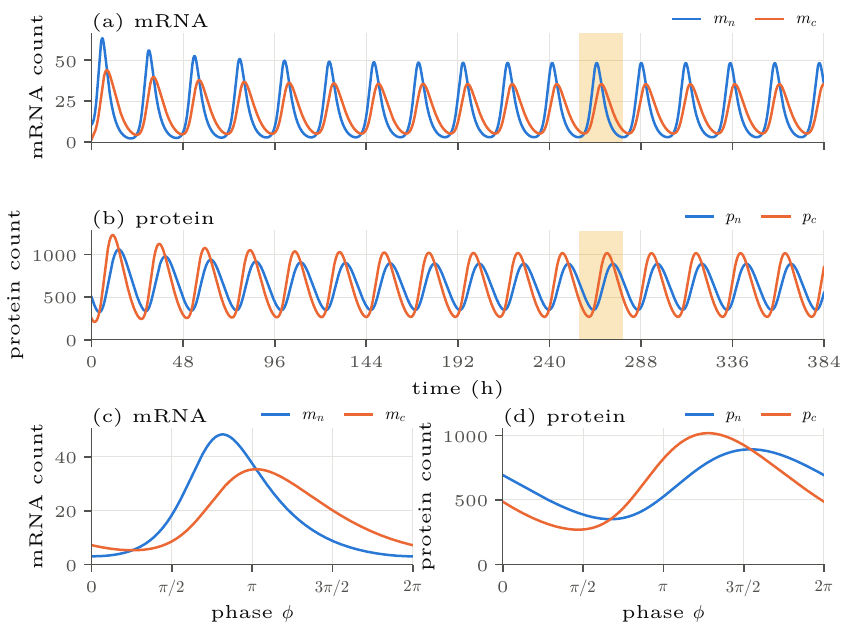}
\caption{Single-cell oscillator with the parameters of \cref{tab:params}. (a, b) Molecule counts against time; the shaded band is the cycle sampled for initial conditions. (c, d) The sampled limit cycle $\mathbf{S}(\phi)$ against phase.}
\label{fig:single}
\end{figure}

\Cref{fig:single} shows the single-cell solution. All four species settle on the limit cycle within about five cycles. On it $m_n$ ranges over $[3.0,\,48.4]$ and $p_n$ over $[351,\,892]$; the Floquet exponents of the cycle are $0$, $-0.0207$, $-0.550$ and $-0.571\ph$, so perturbations of the amplitude relax in about 48 hours.

\subsection{Five cells}

\begin{figure}[!htbp]
\centering
\includegraphics[width=\textwidth]{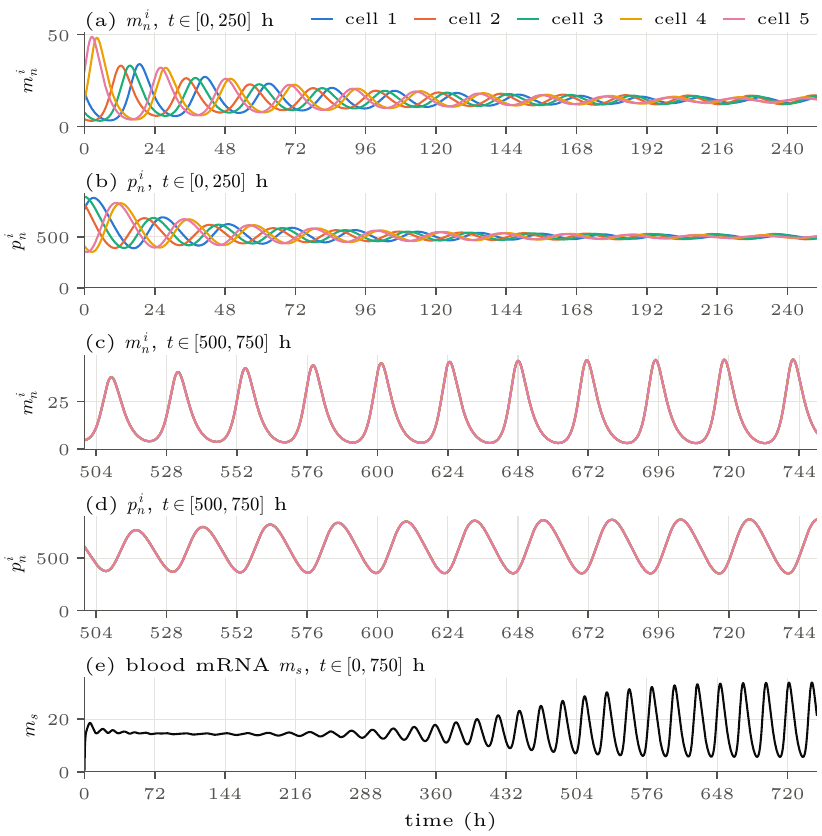}
\caption{Five cells, $\epsilon = 2/22\ph$. (a, b) The phases are spread over the first 250 hours; (c, d) they coincide after 500 hours. (e) The amplitude of the blood level grows as the cells synchronize.}
\label{fig:5cell}
\end{figure}

For five cells and $\epsilon = 2/22\ph$ the phases approach each other slowly (\cref{fig:5cell}): $\spn$ falls from $222$ at $t=0$ to $9.6$ at $250$\,h and $1.2$ at $500$\,h. While the cells are out of phase their contributions to the blood pool nearly cancel, and $m_s$ varies only between $14.3$ and $15.8$ during hours 24 to 72; once the cells are in phase it swings between $5.6$ and $33.9$.

\subsection{One hundred cells}

\begin{figure}[!htbp]
\centering
\includegraphics[width=\textwidth]{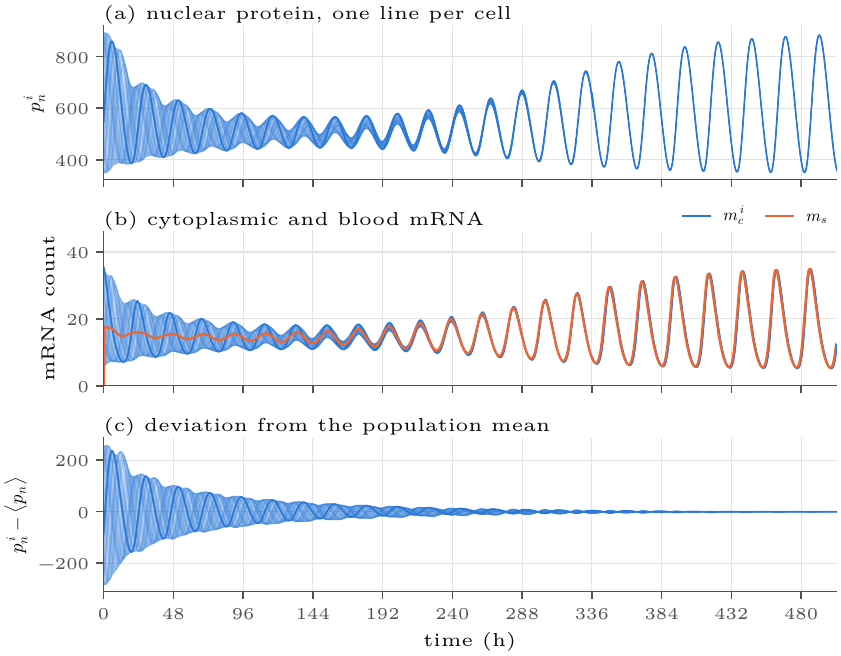}
\caption{One hundred cells, $\epsilon=2/22\ph$, one line per cell. (a) Nuclear protein. (b) Cytoplasmic mRNA of every cell with the blood level $m_s$ drawn over it. (c) Deviation of each cell from the population mean.}
\label{fig:100cell}
\end{figure}

\begin{figure}[!htbp]
\centering
\includegraphics[width=\textwidth]{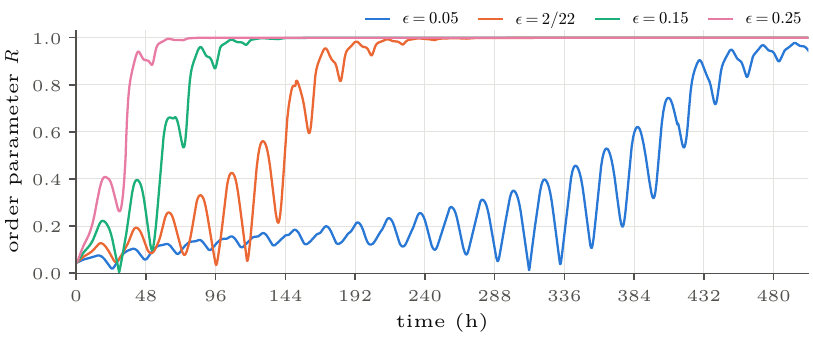}
\caption{Kuramoto order parameter for $100$ cells and four exchange rates, from the same initial phases.}
\label{fig:order}
\end{figure}

\Cref{fig:100cell} shows $100$ cells at the same exchange rate. The spread $\spn$ falls from $195$ to $14$ after $200$\,h and to $0.14$ after $500$\,h. The population mean behaves differently from any single cell: the amplitude of $\E{p_n}$ first collapses, to a swing of only $494$ to $533$ during hours 48 to 96, and then grows back to the single-cell swing of $352$ to $882$ as the cells synchronize. The order parameter rises from $0.04$ to $0.9$ in $167$\,h at this exchange rate, in $428$\,h for $\epsilon=0.05\ph$ and in $40$\,h for $\epsilon=0.25\ph$ (\cref{fig:order}).

\section{Reduction to the population mean}\label{sec:mean}

\Cref{eq:multi} has $4N_c+1$ variables. Averaging over cells gives equations for the means,
\begin{equation}
\begin{aligned}
\dv{\E{m_n}}{t} &= \E{h(p_n)} - \gamma_m \E{m_n}, &\qquad \dv{\E{p_c}}{t} &= \beta \E{m_c} - \gamma_p \E{p_c},\\
\dv{\E{m_c}}{t} &= \gamma_m\E{m_n}-\mu_c\E{m_c}+\mu_s m_s, & \dv{\E{p_n}}{t} &= \gamma_p \E{p_c} -\delta_p \E{p_n},\\
\dv{m_s}{t} &= \mathrlap{N_c\mu_c\E{m_c}-N_c(\mu_s+\delta_m)\, m_s ,}
\end{aligned}
\label{eq:mean-raw}
\end{equation}
which are closed except for the non-linear term $\E{h(p_n)}$.

\subsection{The blood compartment}\label{sec:ms}

\begin{figure}[!htbp]
\centering
\includegraphics[width=\textwidth]{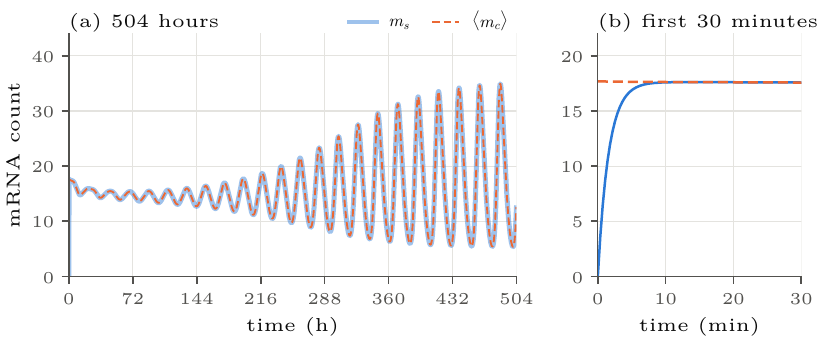}
\caption{Blood level $m_s$ and population mean $\E{m_c}$ for $100$ cells, $\epsilon=2/22\ph$. (a) The two coincide over the whole run. (b) Starting from an empty bloodstream, $m_s$ reaches $\E{m_c}$ within minutes.}
\label{fig:msconv}
\end{figure}

With \cref{eq:mu} the last equation of \eqref{eq:mean-raw} is
\begin{equation}
\dv{m_s}{t}=N_c(\epsilon+\delta_m)\big(\E{m_c}-m_s\big).
\label{eq:ms}
\end{equation}
The blood level relaxes toward $\E{m_c}$ at rate $N_c(\epsilon+\delta_m)$, which for $N_c\gg1$ is much faster than the rate $\sim\nu$ at which $\E{m_c}$ changes. After an initial transient of duration $1/[N_c(\epsilon+\delta_m)]$,
\begin{equation}
m_s=\E{m_c}+O\!\left(\frac{1}{N_c}\right).
\end{equation}
For $N_c=100$ and $\epsilon=2/22\ph$ the relaxation rate is $37.7\ph$. In the simulation the gap falls below $1\%$ of its initial value after $7$ minutes and stays below $1\%$ of $\E{m_c}$ from the first hour on (\cref{fig:msconv}).

Substituting $m_s=\E{m_c}$, the exchange terms cancel in the mean, $-\mu_c\E{m_c}+\mu_s\E{m_c}=-\delta_m\E{m_c}$, while each cell is pulled toward the mean:
\begin{equation}
\dv{m_c^i}{t} = \gamma_m m_n^i-\delta_m m_c^i+\epsilon\big(\E{m_c}-m_c^i\big).
\label{eq:cell-driven}
\end{equation}
The term $\epsilon\,(\E{m_c}-m_c^i)$ is what drives synchronization.

\subsection{Variance correction}

Expanding $h(p_n^i)$ about $\E{p_n}$, the first-order term averages to zero and
\begin{equation}
\begin{aligned}
\E{h(p_n)}&=h\big(\E{p_n}\big)+\tfrac12 h''\big(\E{p_n}\big)\,\spn^2+\dots\\
&=\alpha\hill{\E{p_n}}^r\left[1+\frac{r(r+1)}{2K^2V^2}\hill{\E{p_n}}^{2}\spn^2\right]+\dots
\end{aligned}
\label{eq:taylor}
\end{equation}
To this order the means obey the single-cell equations with one additional term,
\begin{equation}
\begin{aligned}
\dv{\E{m_n}}{t} &= \alpha\hill{\E{p_n}}^r\left[1+\frac{r(r+1)}{2K^2V^2}\hill{\E{p_n}}^{2}\spn^2\right] - \gamma_m \E{m_n},\\
\dv{\E{m_c}}{t} &= \gamma_m\E{m_n}-\delta_m\E{m_c},\qquad
\dv{\E{p_c}}{t} = \beta \E{m_c} - \gamma_p \E{p_c},\qquad
\dv{\E{p_n}}{t} = \gamma_p \E{p_c} -\delta_p \E{p_n}.
\end{aligned}
\label{eq:mean}
\end{equation}
The spread $\spn$ enters as an extra dynamic variable. When $\spn\to0$ the correction vanishes and \eqref{eq:mean} is exactly the single-cell system \eqref{eq:single}: a synchronized population behaves as one cell. Because $h$ is convex, the correction is positive: a desynchronized population transcribes more than a single cell at the mean protein level would.

\begin{figure}[!htbp]
\centering
\includegraphics[width=\textwidth]{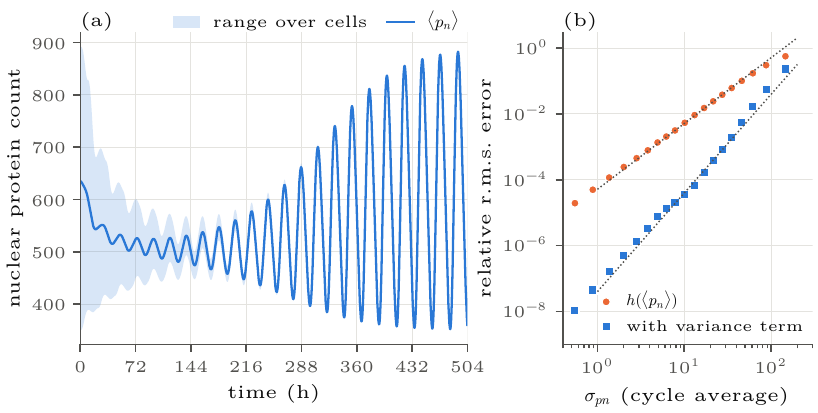}
\caption{One hundred cells, $\epsilon=2/22\ph$. (a) The population mean $\E{p_n}$ against the range covered by the individual cells: the mean rhythm is damped while the cells are incoherent and recovers as they synchronize. (b) Relative r.m.s.\ error over each cycle of the mean transcription rate without and with the variance term of \cref{eq:taylor}, against the cycle-averaged spread. Dotted lines have slopes $2$ and $3$.}
\label{fig:meanfield}
\end{figure}

The expansion parameter is $x=\spn/[V(K+\E{p_n}/V)]$, and the correction in \eqref{eq:taylor} is of relative size $\tfrac12 r(r+1)x^2=15x^2$. At the start of the $100$-cell run $x=0.29$, so the correction is of order one and the series cannot be truncated: the exact mean transcription rate is $4.77\ph$, against $1.31\ph$ without the correction and $3.00\ph$ with it. \Cref{fig:meanfield}(b) compares the two truncations with the exact $\E{h(p_n)}$ cycle by cycle. The error of $h(\E{p_n})$ falls as $\spn^{1.9}$ and that of the corrected expression as $\spn^{3.0}$, as expected from the first neglected term in each case. The corrected expression is within $1\%$ once $\spn<50$ and within $0.1\%$ once $\spn<30$. \Cref{eq:mean} is therefore a description of the approach to synchrony, not of the incoherent state. The collective damping of \cref{fig:meanfield}(a) belongs to the regime where the truncation fails.

\section{Rate of synchronization}\label{sec:rate}

\subsection{Deviations from the mean}

Let $\delta^i=(m_n^i,m_c^i,p_c^i,p_n^i)-\E{(m_n,m_c,p_c,p_n)}$. Subtracting \eqref{eq:mean-raw} from \eqref{eq:multi}, the blood term $\mu_s m_s$ is common to all cells and cancels, so three of the four deviation equations are exactly linear. Only the transcription term needs to be linearized, $h(p_n^i)-\E{h(p_n)}=h'(\E{p_n})\,\delta^i_{pn}+O(|\delta|^2)$. With equal rates $\nu$,
\begin{equation}
\dv{\delta^i}{t}=\big[J(t)-\epsilon P\big]\,\delta^i+O(|\delta|^2),
\qquad
J(t)=\begin{pmatrix}
-\nu & 0 & 0 & h'\big(\E{p_n}(t)\big)\\
\nu & -\nu & 0 & 0\\
0 & \beta & -\nu & 0\\
0 & 0 & \nu & -\nu
\end{pmatrix},
\label{eq:variational}
\end{equation}
where $P=\mathrm{diag}(0,1,0,0)$ selects the $m_c$ component. Near synchrony the mean follows the synchronized orbit, which for $N_c\to\infty$ is the single-cell limit cycle $\mathbf{S}^*(t)$ (\cref{sec:mean}). Then $J(t)$ is the Jacobian of \eqref{eq:single} along $\mathbf{S}^*$ and is $T$-periodic, and \eqref{eq:variational} is a Floquet problem. Every deviation, and hence $\spn$, decays as $e^{-\kappa(\epsilon)t}$ times a $T$-periodic factor, where $-\kappa(\epsilon)$ is the largest Floquet exponent of $J(t)-\epsilon P$.

At $\epsilon=0$ the cells are independent and the largest exponent is $0$: a cell shifted along the limit cycle stays shifted. Coupling makes this neutral direction decay.

\subsection{The first-order rate is \texorpdfstring{$\epsilon/4$}{epsilon/4}}

\begin{proposition}\label{prop:quarter}
Let $\mathbf{S}^*$ be a limit cycle of \eqref{eq:single} with $\gamma_m=\delta_m=\gamma_p=\delta_p=\nu$ whose Floquet exponent $0$ is simple. Then the largest Floquet exponent of $J(t)-\epsilon P$ is
\begin{equation}
-\kappa(\epsilon)=-\frac{\epsilon}{4}+O(\epsilon^2).
\end{equation}
\end{proposition}

\begin{proof}
Consider the more general system $\dd\delta/\dd t=[J(t)-\mathrm{diag}(\epsilon_1,\epsilon_2,\epsilon_3,\epsilon_4)]\,\delta$ and write $\D=\dd/\dd t$. Its four rows read
\[
(\D+\nu+\epsilon_1)\delta_1=h'(p_n^*(t))\,\delta_4,\quad
(\D+\nu+\epsilon_2)\delta_2=\nu\,\delta_1,\quad
(\D+\nu+\epsilon_3)\delta_3=\beta\,\delta_2,\quad
(\D+\nu+\epsilon_4)\delta_4=\nu\,\delta_3 .
\]
Eliminating $\delta_1,\delta_2,\delta_3$ gives a scalar equation for $\delta_4$,
\begin{equation}
\prod_{k=1}^{4}(\D+\nu+\epsilon_k)\;\delta_4=\beta\nu^2\,h'(p_n^*(t))\;\delta_4 .
\label{eq:scalar}
\end{equation}
The change of variables from $\delta$ to $(\delta_4,\D\delta_4,\D^2\delta_4,\D^3\delta_4)$ is linear and invertible. It is also independent of time, because $h'(p_n^*(t))$ first enters at the fourth derivative. It therefore conjugates the monodromy matrix by a constant matrix, so \eqref{eq:scalar} has the same Floquet exponents as the system.

Let $\mu(\epsilon_1,\dots,\epsilon_4)$ denote the largest Floquet exponent, a scalar function on a neighbourhood of the origin in $\mathbb{R}^4$. The coupled model of \cref{eq:variational} is the point $(0,\epsilon,0,0)$, and the uncoupled cell is the origin, where $\mu(0)=0$. We determine the linearization of $\mu$ at the origin. The exponent $0$ is a simple multiplier $\rho=1$ of the monodromy matrix, whose entries depend smoothly on the $\epsilon_k$. Hence $\mu$ is differentiable at the origin, and
\begin{equation}
\mu(\epsilon_1,\dots,\epsilon_4)=\sum_{k=1}^{4}c_k\,\epsilon_k+O(|\epsilon|^2),
\qquad c_k=\left.\frac{\partial\mu}{\partial\epsilon_k}\right|_{\epsilon=0}.
\label{eq:mu-lin}
\end{equation}
Two exact properties of $\mu$ fix the four coefficients.

\emph{Symmetry.} Write $a_k=\nu+\epsilon_k$. The factors $\D+a_k$ have constant coefficients and commute, so
\[
\prod_{k=1}^{4}(\D+a_k)=\D^4+s_1\D^3+s_2\D^2+s_3\D+s_4,
\]
where $s_1=\sum_k a_k$, $s_2=\sum_{j<k}a_ja_k$, $s_3=\sum_{i<j<k}a_ia_ja_k$ and $s_4=a_1a_2a_3a_4$ are invariant under any permutation of the $\epsilon_k$. Equation \eqref{eq:scalar}, and hence $\mu$, is therefore unchanged when the $\epsilon_k$ are permuted. In particular
$\mu(\epsilon,0,0,0)=\mu(0,\epsilon,0,0)=\mu(0,0,\epsilon,0)=\mu(0,0,0,\epsilon)$ for every $\epsilon$. Differentiating at $\epsilon=0$ gives
\begin{equation}
c_1=c_2=c_3=c_4\equiv c .
\label{eq:c-equal}
\end{equation}

\emph{Uniform decay.} On the diagonal $\epsilon_1=\dots=\epsilon_4=\epsilon$ the system is $\dd\delta/\dd t=[J(t)-\epsilon I]\,\delta$. The substitution $\delta=e^{-\epsilon t}\tilde\delta$ gives $\dd\tilde\delta/\dd t=J(t)\,\tilde\delta$, so the monodromy matrix is $e^{-\epsilon T}$ times that at the origin, and every exponent is shifted by exactly $-\epsilon$. Thus $\mu(\epsilon,\epsilon,\epsilon,\epsilon)=-\epsilon$ for every $\epsilon$, and differentiating at $\epsilon=0$ gives
\begin{equation}
\sum_{k=1}^{4}c_k=-1 .
\label{eq:c-sum}
\end{equation}

By \eqref{eq:c-equal} and \eqref{eq:c-sum}, $c=-\tfrac14$, so the linearization is $\mu(\epsilon_1,\dots,\epsilon_4)=-\tfrac14\sum_k\epsilon_k+O(|\epsilon|^2)$. Evaluating it on the coupled model gives $\mu(0,\epsilon,0,0)=-\epsilon/4+O(\epsilon^2)$.
\end{proof}

Both properties hold exactly for every $\epsilon$. The approximation is only in \eqref{eq:mu-lin}: the gradient of $\mu$ at the origin is $-\tfrac14(1,1,1,1)$, but away from the origin $\mu$ is not linear. For the parameters of \cref{tab:params}, the largest exponent with extra decay $0.1\ph$ on any single stage is $-0.0287\ph$ rather than $-0.025\ph$, whereas with $0.1\ph$ on all four stages it is exactly $-0.1\ph$.

\begin{remark}
The result does not depend on $h$, on $\beta$, or on the shape of the limit cycle; it requires only that the four stages have equal first-order rates. In terms of phase reduction~\cite{Kuramoto}, it says that the $m_c$ component carries exactly one quarter of the phase sensitivity of the cycle.
\end{remark}

\begin{remark}\label{rem:fixed}
At the fixed point the right-hand side of \eqref{eq:scalar} is the constant $-a\beta\nu^2$, and the equation reduces to the characteristic polynomial $(\nu+\lambda)^3(\nu+\epsilon+\lambda)+a\beta\nu^2=0$. Differentiating at $\epsilon=0$ gives $4(\nu+\lambda_0)^3\,\delta\lambda+(\nu+\lambda_0)^3\,\delta\epsilon=0$, that is $\delta\lambda=-\epsilon/4$ for every root $\lambda_0$ of \eqref{eq:char}. The same coefficient appears, but the linearization about the fixed point predicts the rate $\epsilon/4-\re\lambda_0$, with $\re\lambda_0=0.041\,\nu>0$ for $G>4$. It applies to cells near the unstable fixed point, not to cells on the limit cycle, where non-linear saturation replaces $\re\lambda_0$ by the neutral exponent $0$.
\end{remark}

\subsection{Numerical results}

\begin{figure}[!htbp]
\centering
\includegraphics[width=\textwidth]{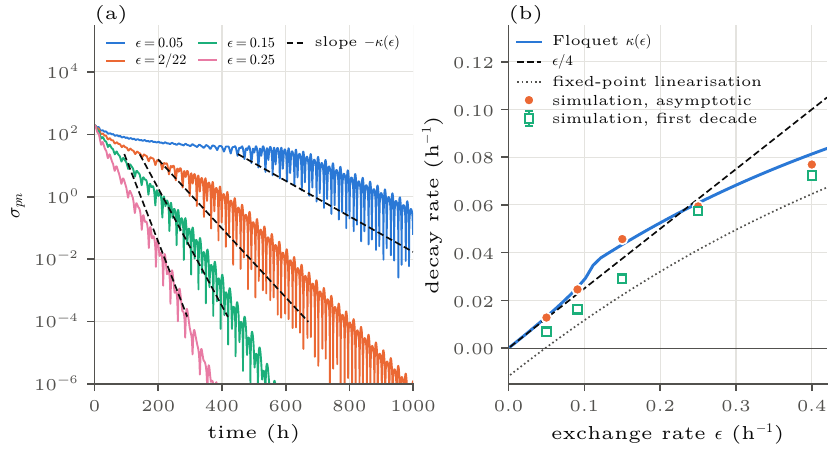}
\caption{Rate of synchronization for $100$ cells. (a) Spread $\spn$ on a logarithmic scale for four exchange rates; dashed segments have the Floquet slope $-\kappa(\epsilon)$. (b) Floquet rate $\kappa(\epsilon)$ on the limit cycle, the first-order result $\epsilon/4$, the prediction of the fixed-point linearization, and rates measured in simulations: asymptotically, and over the first decade of decay from an incoherent start (mean and standard deviation over five sets of initial phases).}
\label{fig:rate}
\end{figure}

\Cref{fig:rate}(b) shows the Floquet rate $\kappa(\epsilon)$ computed for the parameters of \cref{tab:params}. Its slope at the origin is $0.2502$, in agreement with \cref{prop:quarter}. The curve stays close to $\epsilon/4$ over the whole range: $\kappa/\epsilon$ rises to $0.29$ near $\epsilon=0.15\ph$ and falls to $0.20$ at $\epsilon=0.4\ph$. At $\epsilon\approx0.112\ph$ the phase mode meets the amplitude mode of the limit cycle (exponent $-0.0207\ph$ at $\epsilon=0$); beyond this value the two leading multipliers form a complex pair and the decay acquires a slow modulation.

\Cref{tab:rates} compares these rates with simulations of $100$ cells. For $\epsilon\le 2/22\ph$ the measured asymptotic rate agrees with the Floquet exponent of the $100$-cell synchronized orbit to within $0.1\%$. Beyond the mode collision the fitted slopes are within $6\%$; there the slow modulation has a period of tens of cycles, longer than the fitting window. The linearization about the fixed point (\cref{rem:fixed}) underestimates the rate throughout, and predicts no synchronization at all below $\epsilon\approx0.05\ph$.

\begin{table}[!htbp]
\centering\small
\caption{Decay rates of $\spn$ in h$^{-1}$ for $N_c=100$. Floquet exponents are given for the single-cell limit cycle ($N_c\to\infty$) and for the synchronized orbit with $N_c=100$.}
\label{tab:rates}
\begin{tabular}{@{}lcccccc@{}}
\toprule
 & & \multicolumn{2}{c}{\textbf{Floquet}} & \multicolumn{2}{c}{\textbf{Simulation}} & \\
\cmidrule(lr){3-4}\cmidrule(lr){5-6}
$\epsilon$ (h$^{-1}$) & $\epsilon/4$ & $N_c\to\infty$ & $N_c=100$ & asymptotic & first decade & fixed point\\
\midrule
$0.05$ & $0.0125$ & $0.0131$ & $0.0128$ & $0.0128$ & $0.0069\pm0.0015$ & $0.0003$\\
$2/22$ & $0.0227$ & $0.0254$ & $0.0246$ & $0.0246$ & $0.0162\pm0.0010$ & $0.0097$\\
$0.15$ & $0.0375$ & $0.0434$ & $0.0434$ & $0.0457$ & $0.0291\pm0.0002$ & $0.0222$\\
$0.25$ & $0.0625$ & $0.0608$ & $0.0608$ & $0.0594$ & $0.0574\pm0.0006$ & $0.0411$\\
$0.40$ & $0.1000$ & $0.0813$ & $0.0813$ & $0.0769$ & $0.0723\pm0.0016$ & $0.0646$\\
\bottomrule
\end{tabular}
\end{table}

\paragraph{Dependence on population size.} For finite $N_c$ the blood level lags $\E{m_c}$ by about $1/[N_c(\epsilon+\delta_m)]$, which deforms the synchronized orbit and changes the exponent. \Cref{tab:N} shows the effect at $\epsilon=2/22\ph$: the rate rises from $0.0167\ph$ for five cells toward the large-population value $0.0254\ph$, and measured and computed rates agree to within $0.2\%$ at every size.

\begin{table}[!htbp]
\centering\small
\caption{Asymptotic decay rate of $\spn$ against population size, $\epsilon=2/22\ph$.}
\label{tab:N}
\begin{tabular}{@{}lccc@{}}
\toprule
$N_c$ & period of synchronized orbit (h) & Floquet rate (h$^{-1}$) & measured rate (h$^{-1}$)\\
\midrule
$5$ & $23.51$ & $0.01667$ & $0.01668$\\
$20$ & $23.41$ & $0.02213$ & $0.02212$\\
$100$ & $23.36$ & $0.02461$ & $0.02459$\\
$400$ & $23.35$ & $0.02520$ & $0.02516$\\
$\infty$ & $23.34$ & $0.02542$ & --\\
\bottomrule
\end{tabular}
\end{table}

\paragraph{Modulation within a cycle.} In \cref{fig:rate}(a) the spread does not decay smoothly: it dips twice per cycle. For small $\epsilon$ the slowest Floquet mode is a shift along the orbit, $\delta^i\approx-\tau_i\,e^{-\kappa t}\,\dd\mathbf{S}^*/\dd t$ with $\tau_i$ the time offset of cell $i$, so that
\begin{equation}
\spn(t)\approx\sigma_\tau\,e^{-\kappa t}\left|\dv{p_n^*}{t}\right| ,
\end{equation}
which vanishes at the maximum and at the minimum of $p_n$. For $\epsilon=0.05\ph$ the detrended spread $\spn e^{\kappa t}$ follows $|\dd\E{p_n}/\dd t|$ with a correlation of $0.995$ between their logarithms over hours 1000 to 1400.

\paragraph{Approach from the incoherent state.} The exponent $\kappa$ governs the final approach. Starting from uniformly spread phases, the first decade of decay is slower (\cref{tab:rates}): about $0.14\,\epsilon$ for $\epsilon=0.05\ph$ and $0.18\,\epsilon$ for $\epsilon=2/22\ph$. In this phase the mean rhythm is damped (\cref{fig:meanfield}a), so the signal $\E{m_c}$ to which the cells lock is weak. A single exponential fitted to the whole decay therefore gives a rate between the two.

\paragraph{Weak coupling.} Synchrony is not reached from every initial condition. At $\epsilon=0.03\ph$, where $\kappa=0.0077\ph$, four out of four sets of random initial phases led to a state with two groups of cells a few hours apart, with $R\approx0.91$ to $0.94$ and a cycle-averaged $\spn$ of $42$ to $51$ after $4000$ hours, unchanged from its value at $2000$ hours. At $\epsilon=0.05\ph$ all runs synchronized completely. The in-phase state is thus locally stable at rate $\kappa$ but its basin does not contain the incoherent state for all $\epsilon$; we have not located the threshold.

\section{Molecular noise and autocorrelation}\label{sec:noise}

The cells of \cref{sec:results} are deterministic and identical, so once synchronized they stay synchronized. Real cells contain few molecules, here between 3 and 48 nuclear mRNA, and intrinsic noise makes each cell's phase diffuse. We therefore simulate the reactions of \cref{tab:multi-reactions} stochastically for $100$ cells over $2400$\,h (\cref{sec:methods}), and compare independent cells ($\mu_c=\delta_m$, $\mu_s=0$) with cells coupled at $\epsilon=2/22$ and $0.25\ph$. Statistics are taken after $t=800$\,h.

\begin{figure}[!htbp]
\centering
\includegraphics[width=\textwidth]{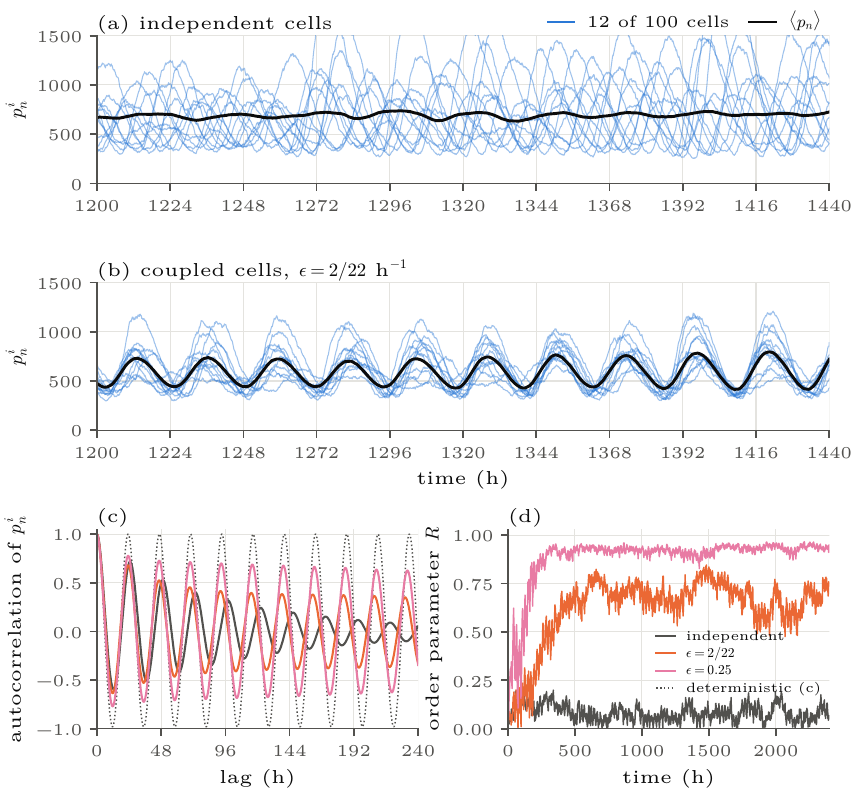}
\caption{Stochastic simulation of $100$ cells. (a, b) Nuclear protein in twelve cells and the population mean over ten days, without and with coupling. (c) Autocorrelation \eqref{eq:acf} of $p_n^i$, averaged over the cells; the dotted curve is the deterministic limit cycle. (d) Order parameter.}
\label{fig:noise}
\end{figure}

\begin{table}[!htbp]
\centering\small
\caption{Stochastic simulations, $N_c=100$, statistics over $t\in[800,2400]$\,h. $C_k$ is the height of the $k$-th peak of the single-cell autocorrelation of $p_n$; s.d.\ of $\E{p_n}$ is the standard deviation in time of the population mean.}
\label{tab:noise}
\begin{tabular}{@{}lcccccc@{}}
\toprule
 & $R$ (mean $\pm$ s.d.) & first peak (h) & $C_1$ & $C_4$ & $C_9$ & s.d.\ of $\E{p_n}$\\
\midrule
Independent cells & $0.08\pm0.04$ & $24.8$ & $0.72$ & $0.32$ & $0.09$ & $25$\\
$\epsilon=2/22\ph$ & $0.68\pm0.07$ & $23.3$ & $0.68$ & $0.42$ & $0.36$ & $113$\\
$\epsilon=0.25\ph$ & $0.93\pm0.02$ & $23.4$ & $0.78$ & $0.70$ & $0.64$ & $172$\\
\bottomrule
\end{tabular}
\end{table}

Independent cells oscillate, but not together (\cref{fig:noise}a): the order parameter stays at the level $\sim N_c^{-1/2}$ of random phases, and the population mean is nearly flat. Their autocorrelation decays from peak to peak and is below $0.1$ at the ninth peak (\cref{tab:noise}); a cell forgets its phase in about ten cycles. Noise also lengthens the mean period, to about $24.3$\,h between successive peaks, against $23.3$\,h for the deterministic cycle.

Coupled cells keep a common rhythm (\cref{fig:noise}b). The autocorrelation of a single cell falls over the first few cycles and then levels off, near $0.36$ for $\epsilon=2/22\ph$ and $0.64$ for $\epsilon=0.25\ph$ at the ninth peak. The plateau is the part of each cell's variance that is locked to the collective rhythm, whose own phase diffuses much more slowly than that of one cell. The order parameter settles at $0.68$ and $0.93$ respectively, and successive peaks are $23.0$ and $23.2$\,h apart, close to the deterministic period. Coupling through the blood therefore both aligns the cells and stabilizes each cell's timing.

\section{Discussion}

Coupling identical feedback oscillators through a shared pool of one species synchronizes them. The rate of the final approach is the Floquet exponent $\kappa(\epsilon)$ of the deviation dynamics, and for a loop of four stages with equal rates it is $\epsilon/4$ to first order, independently of the non-linearity. The argument generalizes directly: for a cyclic chain of $n$ stages with equal first-order rates, exchange of any one species at rate $\epsilon$ gives $\kappa=\epsilon/n+O(\epsilon^2)$.

Several limitations should be kept in mind.
\begin{itemize}[leftmargin=1.4em,itemsep=1pt,topsep=2pt]
\item $\kappa$ is a local rate. The approach from an incoherent population is slower, and at weak coupling the population can remain in a clustered state (\cref{sec:rate}). A global description would require the phase-coupling function of the model, which we have not computed.
\item The variance closure \eqref{eq:mean} is quantitative only for $\spn\lesssim50$; it does not describe the collective damping of the incoherent state.
\item The cells are identical. With a spread of intrinsic periods one expects a threshold in $\epsilon$ below which synchrony is lost~\cite{Kuramoto,Strogatz}; this is not studied here.
\item \Cref{prop:quarter} relies on equal rates for the four stages. With unequal rates the first-order coefficient is no longer $\tfrac14$ and must be computed from the Floquet problem.
\item The blood is treated as a single well-mixed compartment, and the proportionality constants in $\mu_c$, $\mu_s$ and the blood degradation rate are chosen so that synchronized cells reproduce the single-cell model.
\item The noise results are for one population size and one realization per exchange rate.
\end{itemize}

\paragraph{Code availability.} The scripts that produce every figure and number in this paper are included as ancillary files with the arXiv submission:\texttt{sim/model.py}, \texttt{sim/analysis.py} and \texttt{sim/make\_figures.py}.

\appendix
\section{Simulation setups}\label{app:setups}

Unless stated otherwise all runs use the parameters of \cref{tab:params}, $\mu_c=\delta_m+\epsilon$, $\mu_s=\epsilon$, $m_s(0)=0$, and initial phases from \texttt{numpy.random.default\_rng(seed).uniform(0, 2*pi, N)}.

\begin{table}[H]
\centering\small
\caption{All simulation runs.}
\label{tab:all-runs}
\begin{tabular}{@{}llcclcl@{}}
\toprule
\textbf{Run} & \textbf{Type} & $N_c$ & $\epsilon$ (h$^{-1}$) & $t_{\text{end}}$ (h) & seeds & \textbf{Used in}\\
\midrule
Single cell & ODE & $1$ & -- & $504$ & -- & Fig.~\ref{fig:single}\\
Five cells & ODE & $5$ & $2/22$ & $750$ & $7$ & Fig.~\ref{fig:5cell}\\
Hundred cells & ODE & $100$ & $2/22$ & $504$ & $1$ & Figs.~\ref{fig:100cell}, \ref{fig:msconv}, \ref{fig:meanfield}\\
Rate sweep & ODE & $100$ & $0.05$--$0.40$ & $900$--$3800$ & $1$ & Figs.~\ref{fig:order}, \ref{fig:rate}; Table~\ref{tab:rates}\\
Transient rate & ODE & $100$ & $0.05$--$0.40$ & $\le 60/\epsilon+100$ & $1$--$5$ & Fig.~\ref{fig:rate}; Table~\ref{tab:rates}\\
Size sweep & ODE & $5$--$400$ & $2/22$ & $1500$ & $1$ & Table~\ref{tab:N}\\
Weak coupling & ODE & $100$ & $0.03$ & $4000$ & $1$--$4$ & Section~\ref{sec:rate}\\
Noise & tau-leap & $100$ & $0$, $2/22$, $0.25$ & $2400$ & $5$ & Fig.~\ref{fig:noise}; Table~\ref{tab:noise}\\
Step check & tau-leap & $100$ & $2/22$ & $1600$ & $6$ & Section~\ref{sec:methods}\\
\bottomrule
\end{tabular}
\end{table}

\section*{Acknowledgments}
The author would like to thank Charles Peskin for his guidance, encouragement, and valuable discussions throughout the development of this work. Any remaining errors are the sole responsibility of the author.

\end{document}